\documentclass[letterpaper, 10pt, journal, onecolumn, final]{IEEEtran}

\IEEEoverridecommandlockouts                              %

\usepackage{graphicx}

\usepackage{amssymb,amsmath,amsfonts,amsthm}
\usepackage{microtype}
\usepackage{url}
\usepackage[noadjust]{cite}
\usepackage[hidelinks]{hyperref}
\usepackage{xifthen}
\usepackage{xpatch}
\usepackage{cases}
\usepackage[usenames,dvipsnames,table]{xcolor}
\usepackage{tikz}
\usepackage{algorithm}
\usepackage{algpseudocode}
\usepackage{algorithmicx}
\allowdisplaybreaks
\usepackage[normalem]{ulem}
\usepackage{url}
\usepackage{multirow} %
\usepackage{mathtools}
\usepackage{subfig}

\title{\LARGE \bf
A Distributed Bilevel Framework \\for the Macroscopic Optimization of Multi-Agent Systems
}

\author{Riccardo Brumali$^1$, Guido Carnevale$^1$, Sonia Martínez$^2$, Giuseppe Notarstefano$^1$%
\thanks{$^1$ R. Brumali, G. Carnevale, and G. Notarstefano are with Department of Electrical, Electronic and Information Engineering, Alma Mater Studiorum - Universit\`a di Bologna, Bologna, Italy, {\tt\small {name.surname}@unibo.it}}%
\thanks{$^2$ S. Martínez is with Department of Mechanical and Aerospace Engineering, University of California San Diego, La Jolla, CA, USA, {\tt\small soniamd@ucsd.edu}}%
}%

\newcommand{\RB}[1]{{\color{black} #1}}%

\newcommand{\R}{\mathbb{R}}

\newcommand{\N}{\mathbb{N}}

\newcommand{\E}{\mathbb{E}}

\newcommand{\cA}{\mathcal{A}}
\newcommand{\cC}{\mathcal{C}}
\newcommand{\cE}{\mathcal{E}}
\newcommand{\cG}{\mathcal{G}}
\newcommand{\cI}{\mathcal{I}}

\newcommand{\cM}{\mathcal{M}}
\newcommand{\cN}{\mathcal{N}}
\newcommand{\cP}{\mathcal{P}}
\newcommand{\cS}{\mathcal{S}}

\newcommand{\de}{\mathrm{d}}

\newcommand{\T}{^\top}
\newcommand{\inv}{^{-1}}

\newcommand{\cov}{\text{cov}}
\newcommand{\1}{\mathbf{1}}

\newcommand{\argmin}{\mathop{\rm argmin}}

\newcommand{\norm}[1]{\left\lVert#1\right\rVert}

\newcommand\oprocendsymbol{\hbox{$\square$}}
\newcommand\oprocend{\relax\ifmmode\else\unskip\hfill\fi\oprocendsymbol}
\def\eqoprocend{\tag*{$\square$}}

\newtheorem{theorem}{Theorem}

\newtheorem{lemma}{Lemma}
\newtheorem{assumption}{Assumption}
\newtheorem{remark}{Remark}
\newtheorem{example}{Example}

\newcommand{\col}{\text{col}}

\newcommand{\vect}{\text{vec}}
\newcommand{\diag}{\text{diag}}

\newcommand{\tny}[1]{\text{\tiny #1}}

\newcommand{\ud}{^}
\newcommand{\du}{_}

\newcommand{\step}{\gamma}
\newcommand{\stepx}{\gamma_x}
\newcommand{\stepy}{\gamma_y}
\newcommand{\bstepx}{\bar{\step}_{x}}
\newcommand{\bstepxt}{\bar{\step}_{x,1}}
\newcommand{\bstepy}{\bar{\step}_y}

\newcommand{\conv}{\delta}
\newcommand{\bconv}{\bar{\delta}}

\newcommand{\iter}{k}
\newcommand{\iterp}{k+1}

\newcommand{\Proj}{\cP}

\newcommand{\Eig}{\Omega}

\newcommand{\adjKron}{A}

\newcommand{\dimx}{n}
\newcommand{\dimxi}{\dimx\du{i}}
\newcommand{\dimth}{m}

\newcommand{\cCi}{\cC\du{i}}

\newcommand{\Stat}{\Gamma}
\newcommand{\Equi}{\Pi}
\newcommand{\Init}{\Xi}

\newcommand{\uu}{\bar{u}}
\newcommand{\bau}{\bar{v}}
\newcommand{\Sy}{\cS_y}
\newcommand{\Syu}{\cS_y^{\uu}}

\newcommand{\J}{J}

\newcommand{\f}{f}

\newcommand{\g}{g}
\newcommand{\gi}{\g\du{i}}

\newcommand{\p}{p}
\newcommand{\pstar}{\p^{\star}}

\newcommand{\ZZ}{Z}

\newcommand{\w}{w}

\newcommand{\heq}{h_{\tny{eq}}}

\newcommand{\Nc}[1]{\text{N}_{#1}}

\newcommand{\Gfgrad}{G_{\tny{\text{$\!\f$}}}}
\newcommand{\Gthgrad}{G_{\tny{\text{$\!\thh$}}}}
\newcommand{\Gggrad}{G_{\tny{\text{$\!\g$}}}}
\newcommand{\Gghess}{G_{\tny{\text{$\nabla\! \g$}}}}

\newcommand{\LyapF}{U}
\newcommand{\LyapS}{W}
\newcommand{\LyapTot}{V}

\newcommand{\Hess}{H}

\newcommand{\dyn}{\iota}
\newcommand{\dynF}{\dyn_\tny{\text{$B$}}}

\newcommand{\Rmat}{R}

\newcommand{\bi}{b_1}
\newcommand{\bb}{b_2}
\newcommand{\bbb}{b_3}
\newcommand{\bbbb}{b_4}

\newcommand{\ci}{c_1}
\newcommand{\cc}{c_2}
\newcommand{\ccc}{c_3}
\newcommand{\cccc}{c_4}

\newcommand{\Lip}{L}
\newcommand{\Lipf}{\Lip_{\f}}
\newcommand{\Lipmor}{\Lip_{\LyapS}}
\newcommand{\Lipg}{\Lip_{\g}}

\newcommand{\x}{x}
\newcommand{\xt}{\x\ud{\iter}}
\newcommand{\xtp}{\x\ud{\iterp}}
\newcommand{\xii}{\x\du{i}}
\newcommand{\xit}{\xii\ud{\iter}}
\newcommand{\xitp}{\xii\ud{\iterp}}

\newcommand{\bax}{\bar{\x}}

\newcommand{\tx}{\tilde{\x}}
\newcommand{\txi}{\tilde{\x}\du{i}}
\newcommand{\txit}{\txi\ud{\iter}}
\newcommand{\txt}{\tx\ud{\iter}}
\newcommand{\hx}{\hat{\x}}

\newcommand{\thh}{\theta}

\newcommand{\y}{y}
\newcommand{\ystar}{\y^{\star}}
\newcommand{\yt}{\y\ud{\iter}}
\newcommand{\ytp}{\y\ud{\iterp}}
\newcommand{\yi}{\y\du{i}}
\newcommand{\yit}{\yt\du{i}}
\newcommand{\yitp}{\ytp\du{i}}
\newcommand{\yjt}{\yt\du{j}}

\newcommand{\q}{q}
\newcommand{\qt}{\q\ud{\iter}}
\newcommand{\qtp}{\q\ud{\iterp}}
\newcommand{\qi}{\q\du{i}}
\newcommand{\qit}{\qt\du{i}}
\newcommand{\qitp}{\qtp\du{i}}
\newcommand{\qjt}{\qt\du{j}}

\newcommand{\rr}{r}
\newcommand{\rt}{\rr\ud{\iter}}
\newcommand{\rtp}{\rr\ud{\iterp}}
\newcommand{\ri}{\rr\du{i}}
\newcommand{\rit}{\rt\du{i}}
\newcommand{\ritp}{\rtp\du{i}}
\newcommand{\rjt}{\rt\du{j}}

\newcommand{\s}{s}

\newcommand{\sm}{\bar{\s}}
\newcommand{\smt}{\bar{\s}\ud{\iter}}
\newcommand{\smtp}{\bar{\s}\ud{\iterp}}
\newcommand{\sort}{\s_{\perp}}
\newcommand{\sortt}{\s_{\perp}\ud{\iter}}
\newcommand{\sorttp}{\s_{\perp}\ud{\iterp}}

\newcommand{\z}{z}

\newcommand{\zm}{\bar{\z}}
\newcommand{\zmt}{\bar{\z}\ud{\iter}}
\newcommand{\zmtp}{\bar{\z}\ud{\iterp}}
\newcommand{\zort}{\z_{\perp}}
\newcommand{\zortt}{\z_{\perp}\ud{\iter}}
\newcommand{\zorttp}{\z_{\perp}\ud{\iterp}}

\newcommand{\state}{\xi}
\newcommand{\statex}{\state\du{1}}
\newcommand{\statey}{\state\du{2}}
\newcommand{\stater}{\state\du{3}}
\newcommand{\stateq}{\state\du{4}}
\newcommand{\statet}{\state\ud{\iter}}

\newcommand{\stateF}{\zeta}
\newcommand{\stateFt}{\stateF\ud{\iter}}
\newcommand{\stateFtp}{\stateF\ud{\iterp}}
\newcommand{\tstateF}{\tilde{\zeta}}
\newcommand{\tstateFt}{\tilde{\zeta}\ud{\iter}}
\newcommand{\tstateFtp}{\tilde{\zeta}\ud{\iterp}}

\newcommand{\orderLebesgueVal}{6}
\newcommand{\pgraphVal}{0.5}
\newcommand{\wVal}{10^{-4}}
\newcommand{\stepxVal}{10^{-3}}
\newcommand{\stepyVal}{0.9}
\newcommand{\convVal}{0.02}
\newcommand{\dcVal}{0.05}

\newcommand{\orderLebesgue}{\kappa}

\newcommand{\Poly}{P}
\newcommand{\Polyi}{P\du{i}}

\newcommand{\Polyim}{P\du{i-1}}
\newcommand{\Polyimm}{P\du{i-2}}
\newcommand{\bPoly}{\mathbf{P}}

\newcommand{\pgraph}{\tau}

\newcommand{\drift}{d}

\newcommand{\pr}{\alpha}

\newcommand{\fgt}{F_{\stepy}}

\def\Algo/{\RB{Extended Algo}}
\def\algo/{\RB{ALGO}}

\def\Algo/{BILevel Distributed hypergradient for MACRoscopic Optimization}
\def\algo/{BILD-MACRO}

\def\er/{Erd\H{o}s-R\'enyi}

\begin{document}

\maketitle
\thispagestyle{empty}
\pagestyle{empty}

\begin{abstract}
  In this paper, we propose a novel distributed algorithm to optimize the emergent macroscopic behavior of large-scale multi-agent systems via microscopic actions.
  We cast this task as a bilevel optimization problem, where the upper level formalizes the desired macroscopic target behavior through a suitable performance criterion, which is shaped in the lower level by leveraging a compressed aggregate representation estimating the macroscopic state.
   More precisely, the macroscopic state is parametrized by an exponential-family of distributions and constructed from the multi-agent microscopic configuration.
  The proposed algorithm integrates a distributed estimation mechanism, through which each agent reconstructs the macroscopic state locally, with a hypergradient-based update of the microscopic states aimed at improving the collective macroscopic behavior. 
  We prove convergence to the set of stationary points of the bilevel problem via timescale separation arguments. 
  Numerical simulations validate the effectiveness of the proposed method.
\end{abstract}

\section{INTRODUCTION}

Large-scale multi-agent systems are becoming increasingly common in various domains, especially for multi-robot applications~\cite{dorigo2021swarm,testa2025tutorial}.
For these systems, agents interact with each other in order to achieve personal and collective goals, giving rise to complex emergent behaviors at the macroscopic level.
The design of strategies to optimize or impose a desired emergent macroscopic behavior is an important challenge currently studied across several communities~\cite{d2023controlling}. 
A common approach to this problem follows the idea where the microscopic dynamics of the agents are approximated by Partial Differential Equations (PDEs) describing the macroscopic behavior of the system.
Then, the algorithm design is performed at the PDE level and the resulting strategy is implemented microscopically by means of a suitable discretization~\cite{foderaro2014distributed,maffettone2022continuification,maffettone2025leader}.
A direct example of this method is presented in~\cite{di2025continuification}, where a large-scale shepherding problem~\cite{long2020comprehensive} is addressed.
Other applications of this approach include motion planning~\cite{meurer2011finite}, sensor placement~\cite{ghimire2024stein}, and coverage control~\cite{krishnan2022multiscale, todescato2017multi}.
Tightly connected to this strategy, other works cast the problem as an Optimal Transport problem~\cite{bandyopadhyay2014probabilistic,krishnan2018distributed,sinigaglia2021optimal,gao2026banach}.  
In contrast, in~\cite{brumali2025distributed,brumali2025feedback} a different solution is proposed, as agents directly optimize their microscopic state while simultaneously learning the macroscopic behavior.
For large-scale systems, centralized strategies where a single entity controls all the agents are becoming impractical.
Consequently, distributed methods using local computation and neighboring communication are becoming increasingly popular~\cite{giselsson2018large,nedic2018distributed,yang2019survey}.
In this setting, a key aspect concerns the estimation of the macroscopic state which is unknown to individual agents.
We introduce a distributed framework to solve simultaneously both the estimation and optimization tasks leveraging a bilevel optimization formulation~\cite{zhang2024introduction}.
In this context, centralized solution strategies for bilevel optimization problems have been proposed in~\cite{ghadimi2018approximation,pedregosa2016hyperparameter,ji2021bilevel,liu2022bome}, while examples of distributed approaches can be found in~\cite{yousefian2021bilevel,gao2023convergence,niu2026hessian}.
The main contribution of this paper is twofold.
First, we formalize this problem by introducing a distributed bilevel optimization framework tailored for emergent behaviors described by Probability Density Functions (PDFs) that realize a given microscopic configuration.
In our formulation, the lower level estimates a compressed lower-dimensional representation of the macroscopic state from the microscopic configuration, while the upper level enforces a desired macroscopic behavior as a function of the compressed macroscopic state.
By modeling this state as the parameters of an exponential-family distribution, we cast the estimation as a regularized maximum likelihood estimation problem.
Second, we propose \Algo/ (\algo/), a fully distributed method to solve the proposed distributed bilevel formulation.
It integrates a distributed estimation mechanism to locally learn the macroscopic state in each agent, together with a hypergradient-based update on the microscopic states to optimize the emergent behavior.
In turn, we interconnect the estimation and optimization tasks via a consensus mechanism. 
This allows agents to locally reconstruct global quantities from the estimation problem that are required for hypergradient computation.
We highlight that our approach reduces the communication burden with respect to~\cite{brumali2025distributed}.
Indeed, exchanges involve only quantities related to a compressed macroscopic representation instead of a full PDF.

The paper is organized as follows. In Section~\ref{sec:problem_setup} we state the problem setup, in Section~\ref{sec:algo} we present the \algo/ algorithm and its convergence properties, while in Section~\ref{sec:convergence} we theoretically analyze our method. Then, in Section~\ref{sec:simulations} we validate our strategy through numerical simulations. \\%
\paragraph*{Notation}

The symbol $\col(x_1,\dots,x_N)$ denotes the column stacking of the vectors $x_1,\dots,x_N$. $\diag(M_1,\dots,M_N)$ denotes the block diagonal matrix with matrices $M_1,\dots,M_N$ on the diagonal.
$I_m$ denotes the identity matrix in $\R^{m\times m}$.
$1_N$ and $0_N$ denote the vectors of $N$ ones and $N$ zeros, respectively, while %
$\otimes$ is the Kronecker product. 
Dimensions are omitted when they are clear from the context.
We denote the set of positive real numbers as $\R_{+}$. 
Given $v\in \R^n$ and $M \in \R^{n \times n}$, $\norm{v}$ is the Euclidean norm of $v$, while $\norm{M}_F$ is the Frobenius norm of $M$.
Given $f: \R^{n} \times \R^{m} \to \R^n$, $\nabla_1 f(x,y)\in\R^{\dimx}$ and $\nabla_2 f(x,y)\in\R^{\dimth}$ are the gradients of $f$ with respect to its first and second arguments evaluated at $(x,y)$, respectively.
Then $\nabla_{1,2}^2 f(x,y)\in\R^{\dimx \times \dimth}$ and $\nabla_{2}^2 f(x,y)\in\R^{\dimth \times \dimth}$ are the mixed and second-order gradients with respect to the second argument of $f$, respectively.
$\E[\cdot]$ and $\cov[\cdot]$ are the expectation and covariance operators.

\section{Problem Setup and Bilevel Reformulation}\label{sec:problem_setup}
We consider a large-scale multi-agent system composed of $N$ cooperating agents that want to organize themselves to achieve a desired macroscopic behavior. 
Agents can access only their local information and can communicate with each other following the edges of a directed graph $\cG \coloneq \left(\cI, \cE, \cA\right)$.
In particular, we consider $\cI\coloneq \{1,\dots,N\}$ as the set of nodes, i.e., the agents, $\cE \subset \cI \times \cI$ as the set of edges, and $\cA \in \R^{N \times N}$ as the weighted adjacency matrix such that the $(i,j)$-entry $a\du{ij} > 0$ if $(i,j)\in\cE$, otherwise $a\du{ij} = 0$.
We define the set of in-neighbors of agent $i$ as $\cN\du{i} \coloneq \{j \in \cI \mid (j,i) \in \cE\}$.
The next Assumption formalizes the class of graphs we consider.
\begin{assumption} \label{ass:graph}
   The graph $\cG$ is strongly connected and the weighted adjacency matrix $\cA$ is doubly stochastic. \oprocend
\end{assumption}
For each agent $i$, we identify with $\xii \in \cCi \subset \R^{\dimxi}$ its \emph{microscopic} state and we denote with $\x \coloneq \col(\x\du{1},\dots,\x\du{N}) \in \R^\dimx$, where $\dimx \coloneq \sum_{i=1}^{N}\dimxi$, the total microscopic configuration of the system.
We formalize this problem as
\begin{align}\label{eq:general_problem}
   &\min_{\x \in \cC} \int_{\cC} \J(\rho(c,\x)) \de c,
\end{align}
where $\cC \coloneq \cC\du{1} \times \dots \times \cC\du{N}$, $\J: \R_+ \to \R$ is a performance criterion, and $\rho: \cC \times \R^{\dimx} \to \R_+$ is a function describing the macroscopic behavior of the system.

\subsection{Bilevel Reformulation}
Being the function $\rho(\cdot, \x)$ generally unknown, we introduce $\p: \cC \times \R^{\dimth} \to \R_+$, such that
\begin{align*}
   \rho(c,\x) \approx \p(c, \thh(\x)),
\end{align*}
where $\thh: \R^{\dimx} \to \R^{\dimth}$ is a lower-dimensional representation of the microscopic configuration $\x$, while $\p(\cdot, \thh(x))$ is a parametrized family of functions describing the macroscopic behavior of the system.
For example, $\thh(\x)$ can be the result of a model reduction or feature encoding operator, see also Example~\ref{ex:setup}.
We call $\thh(\x)$ the \emph{macroscopic} state, as it is a compressed representation of the macroscopic behavior of the system.
Indeed, we consider cases where $\dimth < \dimx$.

Thus, we reformulate Problem~\eqref{eq:general_problem} by means of a macroscopic cost function $\f : \R^{\dimth} \to \R$ such that 
\begin{align*}
   \f(\thh(\x)) = \int_{\cC} \J(\p(c, \thh(\x))) \de c.   
\end{align*}
Namely, we consider the following optimization problem:
\begin{align}\label{eq:initial_problem}
   \min_{\x\in\cC} \f(\thh(\x)).
\end{align}

To solve Problem~\eqref{eq:initial_problem}, the agents need the knowledge of $\thh(\x)$.
However, as $\thh(x)$ is an emergent large-scale feature of the system, it usually does not have a known expression.
Therefore, it must be estimated from the microscopic configuration $\x$.
To this end, we constrain $\thh(\x)$ to be the solution of an estimation problem and we reformulate problem~\eqref{eq:initial_problem} as the following bilevel optimization problem. Namely, 
\begin{subequations}
   \label{eq:bilevel_problem}
   \begin{align}\label{eq:upper_problem}
       &\min_{\x\in\cC} \quad\f(\thh(\x)) 
       \\
       &\hspace{5pt}\text{s.t.} \quad \ \thh(\x) 
       \in 
       \argmin_{\y \in \R^{\dimth}} \tfrac{1}{N} \textstyle \sum_{i=1}^{N} \gi(\xii, \y), \label{eq:inner_problem}
   \end{align}
\end{subequations}
where $\gi: \R^{\dimxi} \times \R^{\dimth} \to \R$.
For notational convenience, we let $\ell: \R^{\dimx} \to \R$ and $\g : \R^{\dimx} \times \R^{\dimth} \to \R$ 
\begin{align*}
   \ell(\x) \coloneq \f(\thh(\x))
   \quad \text{and} \quad
   \g(\x,\y) \coloneq \tfrac{1}{N}\textstyle\sum_{i=1}^{N} \gi(\xii,\y).
\end{align*}

\begin{example}\label{ex:setup}
   A typical example of this setup is a swarm of robots that want to organize themselves to achieve a certain density in space. 
   The microscopic state of each robot is its position, while the macroscopic state of the system is a parametrization of the robot density in space. 
   In this case, $\thh(\x)$ could be a function associating to the robots spatial configuration a set of parameters of a chosen density model describing the distribution of robots in the environment. \oprocend
\end{example}

In the following assumption, we characterize the property of the set $\cC$ and the function $\f$ that we consider in this paper.
\begin{assumption}\label{ass:upper_cost}
   The set $\cC \subset \R^{\dimx}$ is non-empty, convex, and compact, and 
   $\ell(\x)$ is two-times continuously differentiable.
   \oprocend
\end{assumption}

\subsection{Macroscopic State description}

Following the idea of Example~\ref{ex:setup}, we consider $\thh(\x)$ as a parametrization of the density of the microscopic configuration $\x$.
In detail, we model the agents' density as an exponential-family distribution $\p : \R^{\dimx} \times \R^{\dimth} \to \R_+$. 
Namely, given a parameter vector $\y \in \R^{\dimth}$, we define %
\begin{align}\label{eq:exp_family}
    \p(\x, \y) \coloneq \frac{e^{\y\T \phi(\x)}}{\ZZ(\y)},
 \end{align}
where $\phi : \R^{\dimx} \to \R^{\dimth}$ is a function of sufficient statistics (e.g., a polynomial function of $\x$) and $\ZZ(\y) \coloneq \int_\cC e^{\left(\y\T \phi(c)\right)} \de c$ is a normalization factor.
\begin{remark}\label{rem:p_is_a_pdf}
   Note that, for any fixed $\y \in \R^{\dimth}$, the function $\p(\x, \y)$ is a Probability Density Function (PDF) of $\x \in \cC$. \oprocend
\end{remark}
\begin{remark}\label{rem:phi_is_basis}
   The sufficient statistics $\phi(\x)$ can be interpreted as a set of basis functions that describe the agents' density. 
   The larger its dimension, the more complex densities can be described by the macroscopic state $\thh(x)$. \oprocend
\end{remark}

The next assumption formalizes the class of sufficient statistics we consider in this paper.
\begin{assumption}\label{ass:sufficient_stats}
   The function $\phi(x)$ is two-times continuously differentiable.
   \oprocend
\end{assumption}

In light of Remark~\ref{rem:p_is_a_pdf}, we interpret the estimation problem~\eqref{eq:inner_problem} as a maximum likelihood estimation, where we seek the parameters $\ystar \equiv \thh(\x) \in \R^{\dimth}$ such that $\p(\x, \ystar)$ 
\RB{best fits to} the microscopic configuration $\x$.
Hence, we define $\g(\x, \y)$ as the regularized negative log-likelihood function of $\p(\x, \y)$ given the microscopic configuration $\x$, i.e.,
\begin{align}\label{eq:likelihood_function}
    \g(\x, \y) &\coloneq \tfrac{1}{N}\textstyle\sum_{i=1}^{N} (-\ln (\p(\x\du{i}, \y)) + \frac{\w}{2}\norm{\y}^2)
    \notag\\
   &\ = \tfrac{1}{N}\textstyle\sum_{i=1}^{N} (-\y\T \phi(\x\du{i}) + \ln \ZZ(\y) + \frac{\w}{2}\norm{\y}^2),
\end{align}
where $w>0$ is a regularization parameter.
\section{\algo/ Design and Convergence Properties} \label{sec:algo}
In this section we present \Algo/ (\algo/),
a fully distributed algorithm designed to solve problem~\eqref{eq:bilevel_problem} combining two actions: (i) a learning strategy to estimate the macroscopic state $\thh(\x)$ and (ii) a hypergradient-based update on the microscopic $x$ to solve problem~\eqref{eq:upper_problem}.
In the following, we detail each step and the convergence properties of \algo/, while a formal description of our strategy is presented in Algorithm~\ref{alg:algo_nom}.
It is worth noting that the dimension of the shared quantities in \algo/ does not depend on the agents' number and it scales with the number of sufficient statistics $\dimth$, making our strategy amenable for large-scale systems.
Further, $\dimth$ can be chosen as a trade-off between the accuracy of the macroscopic representation and the computational burden.
\begin{algorithm}%
   \caption{\algo/ -- Agent $i$}
   \label{alg:algo_nom}
   \begin{algorithmic}[0]

      \For{$\iter = 0, 1, 2, \ldots$}

      \State \hspace{-0.4cm}\textbf{Microscopic state update}
      \begin{subequations}
         \label{eq:x_update_algo}
         \begin{align}
            \txit &\!=\! \Proj_{\cC}\!\left[\xit \!-\! \stepx \nabla \phi \! \left(\xit\right) \! \left(\Proj_{\geq \w} \left[\qit\right]\right)\inv \nabla \f\!\left(\yit\right)\right]\! \label{eq:x_update_algo_tilde}
            \\
            \xitp &\!=\! \xit - \conv\left( \txit \!-\! \xit \right) \label{eq:x_update_algo_conv}
         \end{align}
      \end{subequations}

      \State \hspace{-0.4cm}\textbf{Macroscopic learning update}
      \begin{align}\label{eq:y_update_algo}
         \yitp \!=\! \sum_{j\in\cN_i} \!a\du{ij} \yjt \!-\! \stepy \rit
      \end{align}

      \State \hspace{-0.4cm}\textbf{Consensus update}
      \begin{align}
         \ritp &\!=\! \sum_{j\in\cN_i}\! a\du{ij} \rjt \!+\! \nabla_2 \gi(\xitp\!, \yitp) \!-\! \nabla_2 \gi(\xit, \yit)\! \label{eq:r_update_algo}
         \\
         \qitp &\!=\! \sum_{j\in\cN_i}\! a\du{ij} \qjt \!+\! \Hess(\yitp) \!-\! \Hess(\yit) \label{eq:q_update_algo}
      \end{align}

      \EndFor
   \end{algorithmic}
\end{algorithm}

\subsection{Learning the Macroscopic State}
By plugging~\eqref{eq:likelihood_function} into problem~\eqref{eq:inner_problem}, the estimation problem can be rewritten as a distributed consensus optimization %
\begin{align}\label{eq:Learning_problem}
   \min_{\y \in \R^{\dimth}} \tfrac{1}{N} \textstyle\sum_{i=1}^{N} (-\y\T \phi(\x\du{i}) + \ln \ZZ(\y) + \frac{\w}{2} \norm{\y}^2 ).
\end{align}
Let $\gi(\xii, \y) \coloneq -\y\T \phi(\x\du{i}) + \ln \ZZ(\y) + \frac{\w}{2} \norm{\y}^2$.
We solve Problem~\eqref{eq:Learning_problem} in a distributed way \RB{inspired by} a gradient tracking strategy~\cite{nedic2017achieving} where each agent $i$ maintains a local copy $\yi\in\R^{\dimth}$ of the solution estimate $\y$ to problem~\eqref{eq:Learning_problem} and a local proxy $\ri \in \R^{\dimth}$ of the global quantity $\nabla_2 \g(\x, \y)/N$.
Then, at each iteration, agent $i$ updates its own $\yi$ by performing a consensus step with its neighbors and an approximated gradient descent step, i.e., Equation~\eqref{eq:y_update_algo},
where $\stepy > 0$ is a step-size.
Then, the estimate $\ri$ is updated following a dynamic consensus mechanism (see Equation~\eqref{eq:r_update_algo}). 

\subsection{Microscopic State Optimization}
In a centralized setup, an approach to solve problem~\eqref{eq:upper_problem} would be to perform a parallel implementation of a projected gradient descent in which each agent $i$ updates its own $\xii$ as
\begin{subequations}
   \label{eq:centralized_opt_step}
   \begin{align}
      \xitp &= \Proj_{\cC}\left[\xit - \stepx [\nabla\ell(\x)]_i\right] 
      \\
      &= \Proj_{\cC}\left[\xit - \stepx [\nabla\thh(\x)]_i \nabla \f(y)\big|_{\y = \thh(\x)} \right],
   \end{align}
\end{subequations}
where $\stepx > 0$ is a step-size, $\Proj_{\cC}[\cdot]$ is the projection operator onto the set $\cC$, and 
$[\cdot]_i$
is an operator extracting the $i$-th row. 
Note that~\eqref{eq:centralized_opt_step} cannot be implemented because of the lack of analytical knowledge of $\nabla\thh(\x)$.
We overcome this issue following a hypergradient-based approach~\cite[Proposition 1]{ji2021bilevel}.
Specifically, we obtain a local representation of $\nabla\thh(\x)$ around the given point $\x$ by applying the implicit function theorem to the first order optimality condition of Problem~\eqref{eq:inner_problem}, i.e., $\nabla_2 \g(\x, \thh(\x)) = 0$.
It follows that
\begin{align}\label{eq:dthh_dx}
   \nabla\thh(\x) &= -\nabla_2^2 \g(\x,\thh(\x))^{-1} \nabla_{1,2}^2 \g(\x,\thh(\x))\T
   \notag\\
   &= -\Hess(\thh(\x))^{-1}\nabla \phi(\x)\T,
\end{align}
where we consider $\g(\x,\thh(\x))$ as in~\eqref{eq:likelihood_function} and $\Hess(\thh(\x)) \coloneq \nabla_2^2 \g(\x,\thh(\x))$.
However, in our distributed setup, each agent $i$ cannot access the global quantity $\Hess(\thh(\x))$ since it depends on the microscopic $\xii$ of all agents.
We then equip each agent with a local proxy $\qi \in \R^{\dimth \times \dimth}$ of $\Hess(\thh(\x)) / N$ which is updated following a dynamic consensus mechanism, i.e., Equation~\eqref{eq:q_update_algo}.
By combining~\eqref{eq:centralized_opt_step} and~\eqref{eq:dthh_dx} with the local proxy $\qi$, we obtain Equation~\eqref{eq:x_update_algo_tilde},
where $\Proj_{\geq \w} \left[\cdot\right]$ is the projection operator onto the set of symmetric matrices with eigenvalues greater than $\w>0$.
In detail, for a symmetric matrix $\q\in \R^{\dimth \times \dimth}$, we define $\Proj_{\geq \w} \left[\q\right]$ as
\begin{align*}
   \Proj_{\geq \w} \left[\q\right] \coloneq \Eig\max(\w I_\dimth, \Lambda)\Eig\T = \argmin_{B \geq \w I} \norm{A - B}_F^2,
\end{align*}
where $\Eig, \Lambda \in \R^{\dimth \times \dimth}$ are the matrix of eigenvectors and the diagonal matrix of eigenvalues of $\q$, respectively, while the $\max$ operator is applied element-wise (see~\cite[Thm.~2.1]{higham1988computing}).
Then we derive the final update rule~\eqref{eq:x_update_algo_conv} for $\xii$, by performing a convex combination with the previous iterate $\xit$ with parameter $\conv \in (0,1)$.
In particular, as we will show in the analysis, the convex combination is introduced to impose a proper timescale separation between the learning and optimization steps.

\begin{remark}\label{rem:Inverse_Computation}
   The update rule~\eqref{eq:x_update_algo} requires the computation of $\left(\Proj_{\geq \w} \left[\qit\right]\right)\inv$ at each $k>0$, which can be computationally expensive for large $\dimth$.
   A possibility to avoid the explicit computation of the inverse is to solve the linear system $A z = \nabla \f\left(\yit\right)$ for $z \in \R^{\dimth}$ and $A = \Proj_{\geq \w} \left[\qit\right]$. Then $z$ can be used in place of $\left(\Proj_{\geq \w} \left[\qit\right]\right)\inv \nabla \f\left(\yit\right)$ (see~\cite{hestenes1952methods,choi2011minres}). 
   \oprocend
\end{remark}

\subsection{Convergence Properties}
In this section we present the convergence properties of \algo/.
For notational convenience, we rewrite it as %
\begin{subequations}
   \label{eq:compact_sys}   
   \begin{align}
      \!\!\!\xtp &= \xt \!-\! \conv \!\left(\xt \!-\! \Proj_{\cC}\!\!\left[\xt \!-\! \stepx \Gthgrad(\xt\!, \yt\!, \qt) \Gfgrad(\yt)\right]\right)\!\!
      \\
      \!\!\!\ytp &= \adjKron \yt - \stepy \rt
      \\
      \!\!\!\rtp &= \adjKron \rt + \Gggrad(\xtp, \ytp) - \Gggrad(\xt, \yt)
      \\
      \!\!\qtp &= \adjKron \qt + \Gghess(\ytp) - \Gghess(\yt),
   \end{align}
\end{subequations}
where we define the compact matrix $\adjKron \coloneq \cA \otimes I_\dimth$ and the stacked vectors
$\yt \!\coloneq\! \col(\yt\du{1}, \dots , \yt\du{N})$, $\rt \coloneq \col(\rt\du{1}, \dots , \rt\du{N})$, and $\qt \coloneq \col(\qt\du{1}, \dots , \qt\du{N})$.
Further, we consider the stacked operators
$\Gfgrad(\y) \coloneq \col(\nabla \f(\y\du{1}), \dots , \nabla \f(\y\du{N}))$,
$\Gggrad(\x, \y) \!\coloneq\! \col(\nabla_2 \g\du{1}(\x\du{1}, \y\du{1}), \dots , \nabla_2 \g\du{N}(\x\du{N}, \y\du{N}))$,
$\Gthgrad(\x, \y, \q) \!\coloneq\! \diag(-\nabla \phi(\x\du{1}) (\Proj_{\geq \w}[\q\du{1}])\inv, \dots , -\nabla \phi(\x\du{N}) (\Proj_{\geq \w}[\q\du{N}])\inv)$,
and $\Gghess(\y) \!\coloneq\! \col(\Hess(\y\du{1}), \dots , \Hess(\y\du{N}))$.

\RB{Let $\state \coloneq (\statex, \statey, \stater, \stateq) \coloneq (\x, \y, \rr, \q)$ and $\1 \coloneq 1_N \otimes I_\dimth$. 
Then, the set of stationary points of problem~\eqref{eq:bilevel_problem} reads as
\begin{align*}
   \Stat \!\coloneq\! \left\{(\statex, \statey) \!\in\! \cC \times \R^{N\dimth}\!\mid\! \nabla_2 \g(\statex,\statey) \!=\! 0,
   \ \nabla \ell(\statex)\T  \left(\statex \!-\! \statex'\right) \geq 0\   \forall \statex' \in \cC 
   \right\},
\end{align*}
while the set of equilibria of system~\eqref{eq:compact_sys} is defined as
\begin{align*}
   \Equi \!\coloneq\! \left\{\state \!\in\! \Stat \times \R^{N\dimth} \!\times\! \R^{N\dimth \times \dimth} \!\mid 
   \stater \!=\! \tfrac{\1\1\T}{N}\Gggrad(\statex, \1\thh(\statex)),
   \
   \stateq \!=\! \tfrac{\1\1\T}{N}\Hess(\1\thh(\statex)) \!\!\right\}.
\end{align*}
Further, we denote the set of initial conditions as 
\begin{align*}
   \Init \!\coloneq\! \{ \state \!\in\! \cC \times \R^{2N\dimth} \!\times\! \R^{N\dimth \times \dimth} \!\mid\! \ &\stater \!=\! \Gggrad(\statex, \statey),
   \stateq \!=\! \Hess(\statey)\}.
\end{align*}

We state the convergence properties of \algo/.
\begin{theorem}\label{thm:convergence}
   Let Assumptions~\ref{ass:graph},~\ref{ass:upper_cost}, and~\ref{ass:sufficient_stats} hold true. 
   Then for all $\state\ud{0} \in \Init$, %
   there exist $\bstepx, \bstepy, \bconv > 0$ such that for all $\stepx \in (0, \bstepx)$, $\stepy \in (0, \bstepy)$, and $\conv \in (0, \bconv)$, the trajectories of system~\eqref{eq:compact_sys} satisfy
   \begin{align*}
      \lim_{\iter \to \infty} \inf_{\state \in \Equi} \norm{\statet - \state} = 0.
      \eqoprocend
   \end{align*}
\end{theorem}}
\noindent
The proof of Theorem~\ref{thm:convergence} is provided in Section~\ref{sec:proof_thm}.

\section{Convergence Analysis}\label{sec:convergence}

In this section, we theoretically prove Theorem~\ref{thm:convergence}.
Assumptions~\ref{ass:graph}-\ref{ass:sufficient_stats} hold true through the entire section.

\subsection{Preparatory Results}
We start with a preliminary lemma enforcing regularity properties of the cost function $\g(\x, \y)$ (cf.~\eqref{eq:likelihood_function}).
\begin{lemma}\label{lem:regularity_g}
   For all compact set $\Sy \subset \R^{\dimth}$, $\g(\x, \y)$ is $\w$-strongly convex in $\y$ uniformly in $\x \in \cC$ and there exists $\Lipg\!>\!0$ such that $\nabla_2 \g(\x, \y)$, $\Hess(\y)$, and $\nabla \phi(\x)$ are $\Lipg$-Lipschitz continuous and bounded for all $(\x, \y) \in \cC \times \Sy$. 
   \oprocend
\end{lemma}
\noindent
The proof of Lemma~\ref{lem:regularity_g} is provided in Appendix~\ref{sec:proof_regularity_g}.
Then, we introduce the coordinates $\sm \in \R^{\dimth}$, $\zm \in \R^{\dimth\times\dimth}$ and $\sort \in \R^{(N-1)\dimth}$, $\zort \in \R^{(N-1)\dimth\times\dimth}$ defined by 
\begin{align}\label{eq:change}
	\begin{bmatrix}
      \s\\
      \z
   \end{bmatrix} 
   \longmapsto 
   \begin{bmatrix}
      \sm\\
      \sort\\
		\zm\\
		\zort
	\end{bmatrix} \coloneq \begin{bmatrix}
		\tfrac{\1\T}{N} & 0 \\
		\Rmat\T & 0 \\
      0 & \tfrac{\1\T}{N}\\
      0 & \Rmat\T
	\end{bmatrix}
   \begin{bmatrix}
      \rr - \Gggrad(\x, \y) \\
      \q - \Gghess(\y)
   \end{bmatrix}, 
\end{align}
with $\Rmat \in \R^{N\dimth \times (N-1)\dimth}$ such that $\Rmat\T \1 = 0$ and $\Rmat\T \Rmat = I$. 
In light of Assumption~\ref{ass:graph}, it holds $\1\T \adjKron = \1\T$ and $\adjKron\1 = \1$.
Hence, by observing~\eqref{eq:change}, we obtain %
\begin{align*}
   \smtp &= \smt, \quad \zmtp = \zmt,
\end{align*}
leading to $\smt = \sm\ud{0}$ and $\zmt = \zm\ud{0}$ for all $\iter\in\N$. 
Further, the initialization imposed in Theorem~\ref{thm:convergence} implies $\sm\ud{0} = \zm\ud{0} = 0$.
Thus, we neglect $(\smt,\zmt)$ and the resulting system reads as
\begin{subequations}
   \label{eq:compact_sys_causal_no_mean}  
   \begin{align}
      \xtp &\!=\! \xt \!-\! \conv \!\left(\xt \!-\! \Proj_{\cC}\!\!\left[\xt
      - \stepx \Gthgrad(\xt\!, \yt\!, \Rmat \zortt + \Gghess(\yt)) \Gfgrad(\yt)\right]\right) \label{eq:x_update_causal}
      \\
      \ytp &= \adjKron \yt - \stepy \left(\Rmat \sortt + \Gggrad(\xt, \yt)\right) \label{eq:y_update_causal}
      \\
      \sorttp &= \Rmat\T\adjKron\Rmat \sortt + \Rmat\T\left(\adjKron - I\right) \Gggrad(\xt, \yt) \label{eq:s_update_causal}
      \\
      \zorttp &= \Rmat\T\adjKron\Rmat\zortt + \Rmat\T\left(\adjKron - I\right) \Gghess(\yt). \label{eq:z_update_causal}
   \end{align} 
\end{subequations}
System~\eqref{eq:compact_sys_causal_no_mean} enjoys two properties.
First, $\conv \in (0,1)$ allows us to arbitrarily slow down subsystem~\eqref{eq:x_update_causal}. %
Second, the equilibria of subsystem~\eqref{eq:y_update_causal}-\eqref{eq:z_update_causal} are parametrized in $\x$ by the function $\heq : \R^n \to \R^{(2 N -1)\dimth} \times \R^{(N-1)\dimth \times \dimth}$ given by
\begin{align*}
   \heq(\x) \!\coloneq\! \col\!\left(\1\thh(\x), -\Rmat\T\!\Gggrad(\x, \1\thh(\x)), -\Rmat\T\!\Gghess(\1\thh(\x))\right)\!.
\end{align*} 
System~\eqref{eq:compact_sys_causal_no_mean} can be interpreted as a two-time scale system, in which we identify a \emph{fast} subsystem given by the $(\y, \sort, \zort)$-dynamics and a \emph{slow} one given by the $\x$-dynamics.
We analyze them separately.
We start with the so-called \emph{boundary layer system} obtained by considering the $(\y, \sort, \zort)$-error dynamics and by treating the slow variable $\x$ as a constant parameter.
Let $\stateF \coloneq \col(\y, \sort, \zort)$ and $\dynF : \R^n \times \R^{(2 N - 1)\dimth} \times \R^{(N-1)\dimth \times \dimth} \to \R^{(2 N - 1)\dimth} \times \R^{(N-1)\dimth \times \dimth}$ be such that
\begin{align*}
   \stateFtp =\dynF(\x, \stateFt). 
\end{align*}

Then, consider the error coordinates $\tstateF \coloneq \stateF - \heq(\x)$. 
The boundary layer system reads as
\begin{align}\label{eq:err_boundary_layer}
   \tstateFtp = \dynF(\x, \tstateFt + \heq(\x)) - \heq(\x).
\end{align}
We ensure global exponential stability of the origin for~\eqref{eq:err_boundary_layer}.
\begin{lemma}\label{lem:lyapunov_boundary_layer}
   There exists $\LyapF : \R^{(2 N - 1)\dimth} \times \R^{(N-1)\dimth \times \dimth} \to \R$ such that, for all $\tstateF, \tstateF' \!\in\! \R^{(2 N - 1)\dimth} \times \R^{(N-1)\dimth \times \dimth}$, there exists $\bstepy>0$ such that, for all $\stepy \in (0, \bstepy)$, it holds
   \begin{subequations}
      \label{eq:lyapunov_boundary_layer}   
      \begin{align}
         &\bi \|\tstateF\|^2 \leq \LyapF(\tstateF) \leq \bb \|\tstateF\|^2
         \label{eq:lyapunov_boundary_layer_1}   
         \\
         &\LyapF(\dynF(\x, \tstateF + \heq(\x))- \heq(\x)) - \LyapF(\tstateF) \leq -\bbb \|\tstateF\|^2
         \label{eq:lyapunov_boundary_layer_2}   
         \\
         &|\LyapF(\tstateF) - \LyapF(\tstateF')| \leq \bbbb \|\tstateF - \tstateF'\|\left(\|\tstateF\| + \|\tstateF'\|\right),
         \label{eq:lyapunov_boundary_layer_3}   
      \end{align}
   \end{subequations}
   for all $\x \in \cC$ and some $\bi,\bb,\bbb,\bbbb > 0$.
   \oprocend
\end{lemma} 
\noindent
The proof of Lemma~\ref{lem:lyapunov_boundary_layer} is provided in Appendix~\ref{sec:proof_lyapunov_boundary_layer}.
Then, we study the so-called \emph{reduced system}, i.e., the one obtained by plugging $\stateFt = \heq(\xt)$ in the $\x$-dynamics, namely
\begin{align}\label{eq:reduced}
   \xtp \!=\! \xt - \conv \left(\xt - \bax(\xt)\right),
\end{align}
in which $\bax(\xt) \coloneq \Proj_{\cC}[\xt - \stepx \nabla \thh(\xt) \nabla \f\left(\y\right)|_{\y = \thh(\xt)}]$.
Notably, system~\eqref{eq:reduced} corresponds to the ideal $\x$-update~\eqref{eq:centralized_opt_step}.
Next, we ensure LaSalle-type convergence for system~\eqref{eq:reduced}.
\begin{lemma}\label{lem:lyapunov_reduced}
   There exist a continuous $\LyapS \!:\! \R^n \!\to\! \R$ and $\bstepx > 0$ such that, for all $\stepx \in (0, \bstepx)$, it holds
   \begin{subequations}
      \label{eq:lyapunov_reduced}
      \begin{align}
         \!\LyapS(\x - \conv (\x - \bax(\x))) - \LyapS(\x) &\leq - \conv \ci \norm{\x \!-\! \bax(\x)}^2
         \label{eq:lyapunov_reduced_decrease}\!\\
         \norm{\nabla \LyapS(\x)} &\leq \cc \norm{\x - \bax(\x)}
         \label{eq:lyapunov_reduced_nabla_bound}\\
         \norm{\nabla \LyapS(\x) - \nabla \LyapS(\x')} &\leq \ccc \norm{\x - \x'},
         \label{eq:lyapunov_reduced_lip}
      \end{align}
   \end{subequations}
   for all $\x,\x' \in \cC$, $\conv \in (0, 1)$, and some $\ci,\cc,\ccc,\cccc > 0$.
   \oprocend
\end{lemma}
\noindent
The proof of Lemma~\ref{lem:lyapunov_reduced} is provided in Appendix~\ref{sec:proof_lyapunov_reduced}.

\subsection{Proof of Theorem~\ref{thm:convergence}}
\label{sec:proof_thm}

The proof consists of studying the two-time-scale system~\eqref{eq:compact_sys_causal_no_mean} by applying a slightly extended version of the time-scale separation result in~\cite[Thm.~3]{carnevale2025unifying}.
In particular, whereas the cited result establishes global properties under global assumptions, our system enjoys only semi-global properties. 
However, the line of argument is identical and, thus, for reasons of space, we will include this slight extension in a future version of the paper.
More specifically, to apply the semi-global version of~\cite[Thm.~3]{carnevale2025unifying}, we must verify the semi-global version of the three global conditions required in~\cite[Thm.~3]{carnevale2025unifying}.
First, the existence of a Lyapunov function proving exponential stability of the origin for the boundary-layer system~\eqref{eq:err_boundary_layer}, which is guaranteed by Lemma~\ref{lem:lyapunov_boundary_layer}.
Second, the existence of a LaSalle-type function for the reduced system~\eqref{eq:reduced}, which is guaranteed by Lemma~\ref{lem:lyapunov_reduced}.
We point out that the radial unboundedness of $\LyapS$ in our case is not required, since the compact set $\cC$ is forward invariant for the slow system~\eqref{eq:x_update_causal}.
Third, Lipschitz continuity of the dynamics and the equilibrium function $\heq$, and some additional interconnection constraints.
These properties follow from Assumptions~\ref{ass:upper_cost},~\ref{ass:sufficient_stats}, Lemma~\ref{lem:regularity_g}, and the conditions~\eqref{eq:lyapunov_reduced_nabla_bound},~\eqref{eq:lyapunov_reduced_lip} in Lemma~\ref{lem:lyapunov_reduced}.
Hence, all the requirements are satisfied and we introduce $\LyapTot: \cC \times \R^{(2N-1)\dimth} \times \R^{(N-1)\dimth \times \dimth} \to \R$ given by
\begin{align}\label{eq:V}
   \LyapTot(\x, \stateF - \heq(\x)) \coloneq \LyapS(\x) + \LyapF(\stateF - \heq(\x)).
\end{align}
Then, we pick $(\x, \stateF)\!\in\! \cC \!\times\! \R^{(2N-1)\dimth} \!\times\! \R^{(N-1)\dimth \times \dimth}$ and define  
\begin{align}
   \bau \coloneq \LyapTot(\x, \stateF - \heq(\x)) - \inf_{\tau \in \cC}W(\tau),\label{eq:c0}
\end{align}
which exists and is finite since $\cC$ is compact and $W$ is continuous (cf. Assumption~\ref{ass:upper_cost} and Lemma~\ref{lem:lyapunov_reduced}).
By using~\eqref{eq:lyapunov_boundary_layer_1}, we note that the level set $\Omega_{\bau} \coloneq \{\stateF \!\in\! \R^{(2N-1)\dimth} \times \R^{(N-1)\dimth \times \dimth} \mid \LyapF(\stateF - \heq(\x)) \leq \bau\}$ is compact.
Then, we use it to invoke the semi-global version of~\cite[Thm.~3]{carnevale2025unifying} which ensures that for all $(\x, \stateF)\in \cC \times \Omega_{\bau}$ there exists $\bconv \in (0, 1)$ such that, for all $\conv \in (0, \bconv)$, the trajectories of~\eqref{eq:compact_sys_causal_no_mean} satisfy
\begin{align}
   \LyapTot(\x\ud{+}, \stateF\ud{+} - \heq(\x\ud{+})) \!-\! \LyapTot(\x, \stateF - \heq(\x)) 
   \leq - \psi (\|\x\!-\! \bax(\x)\|^2 \!+\! \|\stateF - \heq(\x)\|^2),\label{eq:Lyap_tot}
\end{align}
for some $\psi\!>\!0$, where $(\x\ud{+}, \stateF\ud{+})$ is the update along~\eqref{eq:compact_sys_causal_no_mean} obtained from $(\x,\stateF)$.
By using~\eqref{eq:V},~\eqref{eq:c0}, and~\eqref{eq:Lyap_tot}, it holds
\begin{align}
   \LyapF(\stateF\ud{+} \! - \! \heq(\x\ud{+})) 
   \leq  
   \LyapTot(\x,\stateF \! - \! \heq(\x)) \! - \! \inf_{\tau \in \cC}\LyapS(\tau)
   =
   \bau.
   \label{eq:invariance}
\end{align}
The bound~\eqref{eq:invariance} proves that the compact set $\cC \times \Omega_{\bau}$ is forward invariant for system~\eqref{eq:compact_sys_causal_no_mean}.
Then, for all $(\x\ud0,\stateF\ud0) \in \cC \times \R^{(2N-1)\dimth} \times \R^{(N-1)\dimth \times \dimth}$, we consider $\bau_0 \ge 0$ such that $(\x\ud0,\stateF\ud0) \in \cC \times \Omega_{\bau_0}$. 
Since $\cC \times \Omega_{\bau_0}$ is (i) compact by construction and (ii) forward invariant for system~\eqref{eq:compact_sys_causal_no_mean} by~\eqref{eq:Lyap_tot}, we apply the LaSalle's Invariance Principle~\cite[Thm.~13.3]{haddad2008nonlinear} and prove that the trajectories of~\eqref{eq:compact_sys_causal_no_mean} converge to the largest invariant set $\cM \subset \cC \times \Omega_{\bau_0}$ contained in 
   $E\coloneq \{(\x, \stateF) \in \cC \times \Omega_{\bau_0} \mid \x = \bax(\x), \stateF = \heq(\x)\}$.
We conclude the proof by noting that $\cM \equiv E$.

\section{Numerical Simulations} \label{sec:simulations}
\begin{figure}[htpb]
   \centering
   \centering
   \subfloat[]{
      \begin{tabular}{cccc}
         \includegraphics[width=0.23\textwidth]{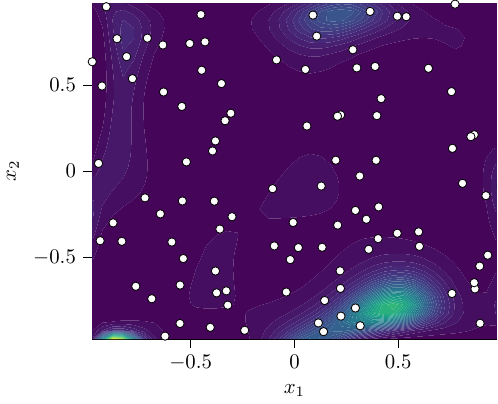} & 
         \includegraphics[width=0.23\textwidth]{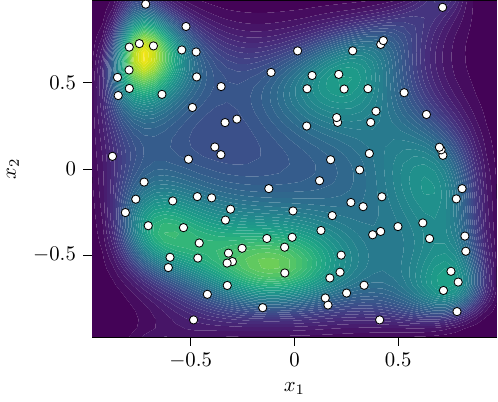} & 
         \includegraphics[width=0.23\textwidth]{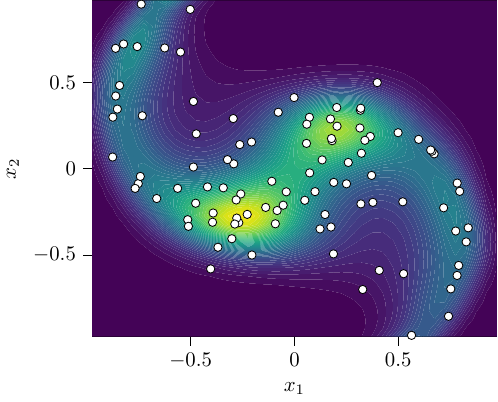} & 
         \includegraphics[width=0.23\textwidth]{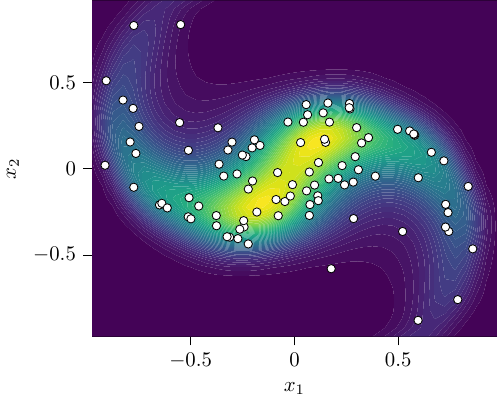} \\
         \scriptsize iteration $\iter=0$ & 
         \scriptsize iteration $\iter=200$ & 
         \scriptsize iteration $\iter=2000$ & 
         \scriptsize iteration $\iter=100000$
      \end{tabular}
   }

   \subfloat[]{
      \begin{tabular}{cccc}
         \includegraphics[width=0.23\textwidth]{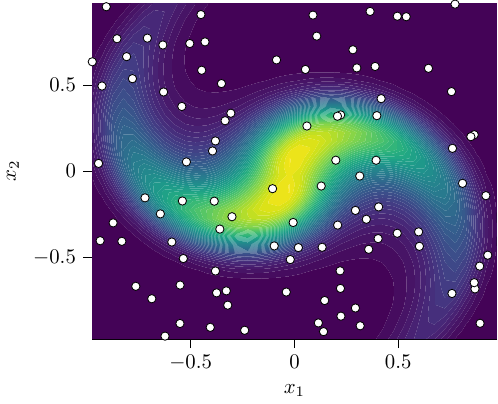} & 
         \includegraphics[width=0.23\textwidth]{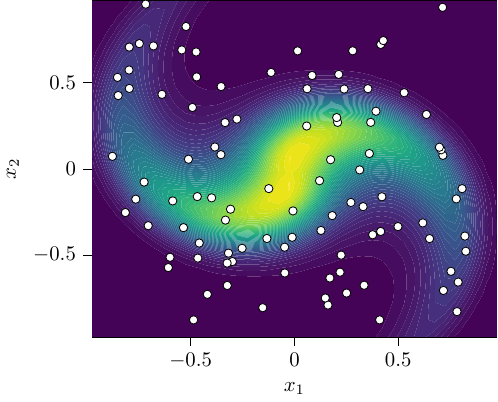} & 
         \includegraphics[width=0.23\textwidth]{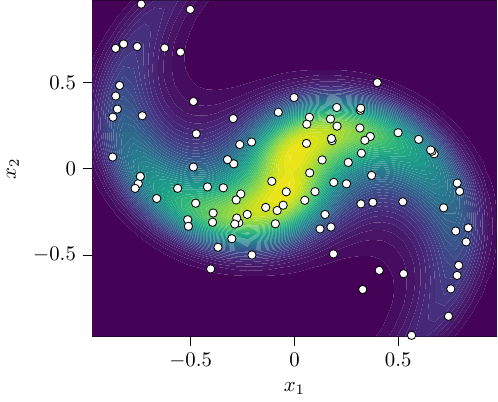} & 
         \includegraphics[width=0.23\textwidth]{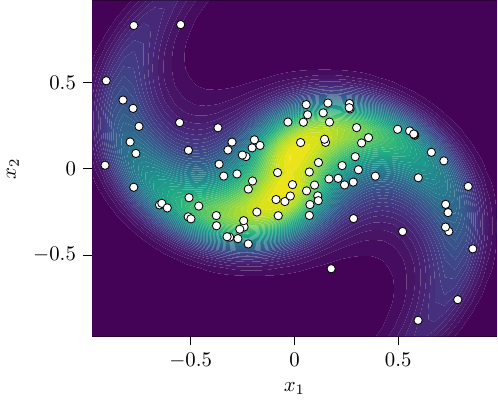} \\
         \scriptsize iteration $\iter=0$ & 
         \scriptsize iteration $\iter=200$ & 
         \scriptsize iteration $\iter=2000$ & 
         \scriptsize iteration $\iter=100000$
      \end{tabular}
      \label{fig:frames_star}
   }
   \caption{Evolution of the agents' states (white dots) at different iterations $\iter$, with (a) the learned density and (b) the target density in background.}
   \label{fig:frames}
\end{figure}

In this section, we test the effectiveness of \algo/ through numerical simulations on a multi-robot setup. 
We consider a swarm of $N$ robots that want to organize themselves to mimic a desired space density.
We define each $\cCi \coloneq [-1, 1] \times [-1, 1]$ and we choose the microscopic state $\xii \in \R^2$ of each robot as its position for all $i \in \cI$.
To describe the density of the robots in space, we consider as sufficient statistics for $\p(\cdot, \thh(x))$ (cf.~\eqref{eq:exp_family}) the Legendre polynomials $\Polyi: \R \to \R$ up to order $\orderLebesgue = \orderLebesgueVal$ (see~\cite{lima2022lecture}), namely 
\begin{align*}
   \Polyi(c) \coloneq  \tfrac{1}{i}\left(\left(2 i - 1\right) c \Polyim - \left(i - 1\right) \Polyimm\right), \quad  i \geq 2,
\end{align*}
with $\Poly\du{0}(c) \!=\! 1$ and $\Poly\du{1}(c) \!=\! c$.
Then, we define $\phi(c) \in \R^{\dimth}$ as
\begin{align*}
   \phi(c) = \vect[\bPoly(c) \otimes \bPoly(c)],  
\end{align*}
with $\bPoly(c) \!\coloneq\! \col(\Poly\du{0}(c), \dots, \Poly\du{\orderLebesgue}(c))$ and $\vect[\cdot]$ as the vectorization operator.
We define $\f$ as the Kullback-Leibler divergence between $\p(\cdot, \thh(\x))$ and a desired $\pstar$ (see Fig.~\ref{fig:frames_star}), i.e.,
\begin{align}\label{eq:Kullback_Leibler}
   \f\left(\thh(\x)\right) = \int_{\cC} \p\left(c, \thh(\x)\right) \ln \left(\tfrac{\p\left(c, \thh(\x)\right)}{\pstar(c)}\right) \de c.
\end{align}
We point out that the evaluation of~\eqref{eq:Kullback_Leibler} is done through numerical integration by discretizing the set $\cC$ in a grid of equally spaced points with discretization step $\Delta c = \dcVal$. 
For all $i \in \cI$, the initial $\xii\ud0$ and $\yi\ud0$ are sampled from a Uniform and Normal distribution, respectively.
For the communication graph $\cG$, we consider an \er/ graph with edge probability $\pgraph = \pgraphVal$.
As for the algorithm parameters we consider $\stepx = \stepxVal$, $\stepy = \stepyVal$, $\w = \wVal$, and $\conv = \convVal$.

\begin{figure}[htpb]
   \centering
   \subfloat[][]{\label{fig:descent_lower}\includegraphics[scale=0.44]{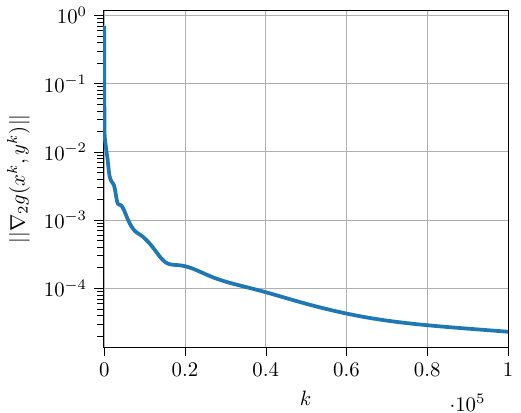}} \ 
   \subfloat[][]{\label{fig:descent_upper}\includegraphics[scale=0.44]{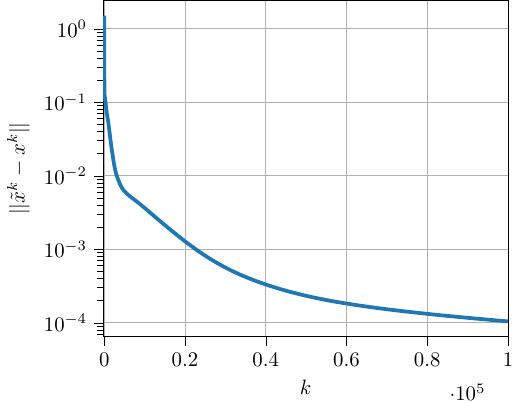}} 
   \caption{Evolution of (a) $\norm{\nabla_2 \g(\xt, \yt)}$ and (b) $\norm{\txt - \xt}$.}
   \label{fig:descent}
\end{figure}

In Fig.~\ref{fig:descent_lower} we report the evolution of $\norm{\nabla_2 \g(\xt, \yt)}$, while Fig.~\ref{fig:descent_upper} shows the evolution of $\norm{\txt - \xt}$, where $\txt \coloneq \col(\txt\du{1}, \dots, \txt\du{N})$.
In particular, both quantities converge to 0 implying the convergence of the sequences $\{\xt\}_{\iter\in\N}$ and $\{\yt\}_{\iter\in\N}$ to a stationary point of problem~\eqref{eq:bilevel_problem}. 

Then, in Fig.~\ref{fig:frames}, we provide snapshots of the robots configuration during the execution of \algo/ showing how the swarm organizes itself and learns the macroscopic density to mimic the desired $\pstar$.
In particular, in Fig.~\ref{fig:frames}(a) we report the evolution of the agents' states (white dots) at different iterations $\iter$, with the learned density $\p(\cdot, \yit)$ in background, while in Fig.~\ref{fig:frames}(b) we show the same agents' states but with the desired $\pstar$ in background.
It is possible to observe that the evolution of the agents' states and the learned $\p(\cdot, \yit)$ change in order to mimic the desired $\pstar$.

\section{Conclusions}
In this paper we proposed \algo/, a fully distributed algorithm to optimize the macroscopic behavior of large-scale multi-agent systems.
We introduced a bilevel optimization formulation and designed \algo/ by combining a distributed consensus strategy for the learning problem, and a distributed hypergradient update on the microscopic state. %
We proved the convergence of \algo/ to stationary points of the problem and we validated the proposed approach through numerical simulations.
\appendix
\subsection{Proof of Lemma~\ref{lem:regularity_g}}\label{sec:proof_regularity_g}
   We start by computing $\nabla_2 \g(\x,\y)$, namely
   \begin{align}\label{eq:grad_g}
      \!\!\!\nabla_2 \g\left(\x,\y\right) \!\coloneq\! \frac{1}{N}\sum_{i=1}^{N} \!\left(\!-\phi(\xii) \!+\! \!\! \int_{\cC}\! \p(c, \y) \phi(c) \de c \!+\! \w \y\right).
   \end{align}

   \noindent
   Then, we evaluate the Hessian $\Hess(\y)$ as 
   \begin{subequations}
      \begin{align*}
         \Hess(\y)
         &=  - \left(\int_{\cC} \p(c, \y) \phi(c) \de c \right)\left(\int_{\cC} \p(c, \y) \phi(c) \de c \right)\T 
		+  \int_{\cC} \p(c, \y) \phi(c)\phi(c)\T \de c + \w I
         \notag\\
         &\stackrel{(a)}{=}  \E\left[\phi(c)\phi(c)\T\right] 
         - \E\left[\phi(c)\right]\E\left[\phi(c)\right]\T + \w I
         \notag\\
         &\stackrel{(b)}{=}  \cov\left[\phi(c)\right] + \w I,
      \end{align*}
   \end{subequations}
   where in $(a)$ we use the definition of the expectation operator while in $(b)$ we use the property $\cov\left[\phi(c)\right] = \E\left[\phi(c)\phi(c)\T\right] - \E\left[\phi(c)\right]\E\left[\phi(c)\right]\T$.
   Since $\cov\left[\phi(c)\right] \!\geq\! 0$ by definition, then $\Hess(\y) \!\geq\! \w I$, implying the $\w$-strong convexity of $\g(\x, \y)$ with respect to $\y$ for all $\x \in \cC$ and $\y \in \Sy$.
   Finally, the Lipschitz continuity and boundedness of $\nabla_2 \g(\x, \y)$, $\Hess(\y)$, and $\nabla\phi(\x)$ follow by their continuous differentiability (cf. Assumption~\ref{ass:sufficient_stats}) and the compactness of $\cC \times \Sy$.

\subsection{Proof of Lemma~\ref{lem:lyapunov_boundary_layer}}\label{sec:proof_lyapunov_boundary_layer}
   The system~\eqref{eq:err_boundary_layer} is a cascade of two subsystems. 
   First, a reformulation of the gradient tracking algorithm expressed in suitable coordinates with an arbitrarily fixed $x \in \cC$, and,
   second, a dynamic consensus algorithm.
   We highlight this structure by rewriting system~\eqref{eq:err_boundary_layer} as 
   \begin{subequations}\label{eq:cascade_system}
      \begin{align}
         \begin{bmatrix}
            \tstateFtp_1 
            \\
            \tstateFtp_2 
         \end{bmatrix} &= \fgt(\x, \tstateFt_1,\tstateFt_2)
         \label{eq:GT} 
         \\
         \tstateFtp_3 &= \Rmat\T A \Rmat \tstateFt_3 + \drift(\x, \tstateFt_1),
         \label{eq:dynamic_consensus}
      \end{align}
   \end{subequations}
   where $\fgt: \cC \times \R^{(2N-1)\dimth} \to \R^{(2N-1)\dimth}$ is the state transition function of the $(\tstateF_1,\tstateF_2)$-dynamics, while $\drift: \cC \times \R^{(N-1)\dimth} \to \R^{(N-1)\dimth \times \dimth}$ captures the coupling between the two subsystems due to the drift of $\tstateFt_1$, namely 
      \begin{align*}
         \fgt(\x,\tstateFt_1,\tstateFt_2) &\!\coloneq\! 
         \begin{bmatrix}
            \adjKron \tstateFt_1 - \stepy \left(\Rmat\tstateFt_2 + \Gggrad(\x,\tstateFt_1)\right)  
            \\
            \Rmat\T \!\adjKron\Rmat \tstateFt_2 + \Rmat\T \!\left(\adjKron - I\right) \Gggrad(\x,\tstateFt_1 + \thh(\x))
         \end{bmatrix}
         \\
         \drift(\x,\tstateFt_1) &\!\coloneq \! \Rmat\T\!\!\left(\adjKron \!-\! I\right) \! \big(\Gghess(\tstateFt_1 \!+\! \1\thh(x)) \!-\! \Gghess(\1\thh(x))\big).
      \end{align*}%
   We analyze System~\eqref{eq:GT} with the same Lyapunov function used in~\cite{notarnicola2023gradient} for a single-problem context, namely 
   \begin{align}\label{eq:lyapunov_eq}
      \LyapF_1(\tstateF_1,\tstateF_2) \coloneq 
      \begin{bmatrix}
         \tstateF_1\\
         \tstateF_2
      \end{bmatrix}\T P_{\tny{GT}}\begin{bmatrix}
         \tstateF_1\\
         \tstateF_2
      \end{bmatrix},
   \end{align}
   where $P_{\tny{GT}} = P_{\tny{GT}}\T \in \R^{(2N-1)\dimth \times (2N-1)\dimth}$ is a positive definite matrix, see~\cite{notarnicola2023gradient} for its specific definition. 
   For any $\uu > 0$ and $\x \in \cC$, 
   let us introduce the compact set $\Syu$ defined as
   \begin{align*}
      \Syu \coloneq \{(\stateF_1,\stateF_2) \in \R^{(2N-1)\dimth} \mid
      \LyapF_1(\stateF_1 - \1\thh(x), \ \stateF_2 + \Rmat\T \Gggrad(\x,\1\thh(x))) \leq \uu \}.
   \end{align*} 
   By Lemma~\ref{lem:regularity_g}, for all $\uu > 0$ and $(\x,\stateF) \in \cC \times \Syu \times \R^{(N-1)\dimth\times\dimth}$, we can ensure that $\g(x,\stateF_1)$ is $\w$-strongly convex in $\stateF_1$ with $\Lipg$-Lipschitz continuous gradients.
   By recalling that $\cG$ is strongly connected and $\cA$ is doubly stochastic (cf. Assumption~\ref{ass:graph}), we can slightly adapt~\cite[Thm.~4.3]{notarnicola2023gradient} to show that for all $\uu>0$ and $(\stateF_1, \stateF_2) \in \Syu$, there exists $\bar{\step}_y > 0$ such that, for all $\stepy \in (0,\bar{\step}_y)$, we can bound the increment $\Delta \LyapF_1(\tstateF_1,\tstateF_2) \coloneq \LyapF_1(\fgt(\x,\tstateF_1,\tstateF_2)) - \LyapF_1(\tstateF_1,\tstateF_2)$ of $\LyapF_1$ along the trajectories of system~\eqref{eq:GT} as 
   \begin{align}
      \Delta \LyapF_1(\tstateF_1,\tstateF_2) \leq - b_{2,1}\norm{\begin{bmatrix}\tstateF_1\\\tstateF_2\end{bmatrix}}^2,
      \label{eq:increment_lyapF_1}
   \end{align}
   for some $\tilde{b}_{2,1} > 0$.
   Then, we study subsystem~\eqref{eq:dynamic_consensus} by introducing $\LyapF_2: \R^{(N-1)\dimth\times\dimth} \to \R$
   \begin{align}
      \LyapF_2(\tstateF_3) \coloneq \tstateF_3\T P\tstateF_3,
   \end{align}
   where $P = P\T \in \R^{(N-1)\dimth \times (N-1)\dimth}$ is the solution to
   \begin{align}
      (\Rmat\T A\Rmat)\T P \Rmat\T A\Rmat - P = -I.
   \end{align}
   Being $\Rmat\T A\Rmat$ Schur (cf. Assumption~\ref{ass:graph}), $P$ is also unique and positive definite. 
   Let $\Lip_\drift>0$ be the Lipschitz continuity constant of $\drift(\x,\tstateF_1)$.
   The increment $\Delta \LyapF_2(\tstateF_3) \coloneq \LyapF_2(\Rmat\T A \Rmat \tstateF_3 + \drift(\x,\tstateF_1)) - \LyapF_2(\tstateF_3)$ of $\LyapF_2$ along the trajectories of system~\eqref{eq:dynamic_consensus} reads as 
   \begin{align}
      \Delta \LyapF_2(\tstateF_3) 
      &\!=\! \tstateF_3\T (\Rmat\T\!\! A\Rmat)\T\! P (\Rmat\T\!\! A\Rmat) \tstateF_3 - \tstateF_3\T P \tstateF_3 
      \notag\\ 
      &\hspace{9pt}
      \!+\! 2 \tstateF_3\T (\Rmat\T\!\! A\Rmat)\T\! P \drift(\x,\tstateF_1) 
      + \drift(\x,\tstateF_1)\T\! P \drift(\x,\tstateF_1)
      \notag\\
      &\stackrel{(a)}{\leq} -\|\tstateF_3\|^2
      + 2 \Lip_\drift\|\tstateF_3\| \|\left(\Rmat\T\!\! A\Rmat\right)\T\!\! P \| \norm{\begin{bmatrix}\tstateF_1 \\ \tstateF_2\end{bmatrix}}
      \notag\\ 
      &\hspace{.4cm}
      + \Lip_\drift^2\norm{P}\norm{\begin{bmatrix}\tstateF_1 \\ \tstateF_2\end{bmatrix}}^2,
      \label{eq:increment_lyapF_2}
   \end{align}
   where in $(a)$ we apply~\eqref{eq:lyapunov_eq} and the Cauchy-Schwarz inequality, 
   we use the Lipschitz continuity of $\drift(\x,\tstateF_1)$, and the fact that $\drift(\x,0) = 0$.
   With these results at hand, we now introduce the Lyapunov function $\LyapF: \R^{(2N-1)\dimth} \times \R^{(N-1)\dimth \times \dimth} \to \R$
   \begin{align}
      \LyapF(\tstateF_1,\tstateF_2,\tstateF_3) \coloneq \LyapF_1(\tstateF_1,\tstateF_2) + \pr \LyapF_2(\tstateF_3),
   \end{align}
   where $\pr > 0$ is a design parameter to be selected.
   Since $\LyapF$ is quadratic and $P_{\text{\tny{GT}}}, P > 0$, the bounds~\eqref{eq:lyapunov_boundary_layer_1} and~\eqref{eq:lyapunov_boundary_layer_3} are trivially verified.
   To verify~\eqref{eq:lyapunov_boundary_layer_2}, let us introduce $t_1 \coloneq \RB{\Lip_\drift^2}\norm{P}$ and $t_2 \coloneq 2 \Lip_\drift\norm{\Rmat\T \adjKron \Rmat P}$. 
   By using~\eqref{eq:increment_lyapF_1} and~\eqref{eq:increment_lyapF_2}, the increment of $\Delta \LyapF(\tstateF_1,\tstateF_2,\tstateF_3) \coloneq \Delta \LyapF_1(\tstateF_1,\tstateF_2) + \pr \Delta \LyapF_2(\tstateF_3)$ along the trajectories of system~\eqref{eq:cascade_system} satisfies
   \begin{align*}
      \Delta \LyapF(\tstateF_1,\tstateF_2,\tstateF_3) \!\leq\! 
      - \!\! \begin{bmatrix}
         \norm{\begin{bmatrix}
         \tstateF_1\\\tstateF_2
         \end{bmatrix}}
         \\
         \|\tstateF_3\|
      \end{bmatrix}\T  
      \!\!\!\underbrace{\begin{bmatrix}
         b_{2,1} \!-\! \alpha t_1& -\alpha t_2
         \\
         -\alpha t_2 &\alpha
      \end{bmatrix}}_{Q(\alpha)}
      \!\!\begin{bmatrix}
         \norm{\begin{bmatrix}
         \tstateF_1\\\tstateF_2
         \end{bmatrix}}
         \\
         \|\tstateF_3\|
      \end{bmatrix}\!\!,
   \end{align*}
   for all $\tstateF \in \R^{(2N-1)\dimth} \times \R^{(N-1)\dimth \times \dimth}$ such that $(\tstateF_1 + \1\thh(x), \tstateF_2 - \Rmat\T \!\Gggrad(\x,\1\thh(x))) \in \Syu$. 
   By selecting $\pr \in (0,\bar{\pr})$ with $\bar{\pr} \coloneq \min\{\frac{b_{2,1}}{t_1}, \frac{b_{2,1}}{t_1 - t_2^2}\}$, the Sylvester criterion ensures the positive definiteness of $Q(\alpha)$. Thus, also~\eqref{eq:lyapunov_boundary_layer_2} is achieved.

\subsection{Proof of Lemma~\ref{lem:lyapunov_reduced}}\label{sec:proof_lyapunov_reduced}
   We start by showing that $\nabla \ell(\x) \coloneq \nabla \thh(\x) \nabla \f(y)\!\!\mid_{y=\thh(\x)}$ is Lipschitz continuous for all $\x \in \cC$.
   In light of the compactness of $\cC$ (cf. Assumption~\ref{ass:upper_cost}), the continuous differentiability property of $\nabla\phi(\x)$ (cf. Assumption~\ref{ass:sufficient_stats}) implies the boundedness and Lipschitz continuity of $\nabla\phi(\x)$ for all $\x \in \cC$.  
   Hence, the Lipschitz continuity of $\thh(\x)$ follows by~\cite[Lemma 2.2]{ghadimi2018approximation}.
   Then, for all $\x \in \cC$, the Lipschitz continuity and boundedness of $\Hess(\thh(\x))$ follow since it is obtained from the composition of continuously differentiable functions.
   Similarly, it can be shown that $\nabla f(\y)\!\!\mid_{y=\thh(\x)}$ is Lipschitz continuous and bounded.
   Thus, the $\Lipf$-Lipschitz continuity of $\nabla \ell(\x)$ follows since $\ell$ is the product of Lipschitz continuous and bounded functions. 
   Finally, we point out that the Lipschitz continuity of $\nabla \ell(\x)$ implies the $\Lipf$-weak convexity of $\ell(\x)$, for all $\x \in \cC$.
   Then, we select $\LyapS: \R^{\dimx} \to \R$ as the Moreau envelope of $\ell(\x)$, i.e.,
   \begin{align*}
      \LyapS(\x) \coloneq \inf_{\y \in \cC} \left(\ell(\y) + \tfrac{1}{2\stepx} \norm{\y - \x}^2\right).
   \end{align*}
   For all $\stepx \!\in\! (0, \bstepxt)$ with $\bstepxt \!\coloneq\! \tfrac{1}{\Lipf}$, $\LyapS(\x)$ has $\Lipmor$-Lipschitz continuous gradient~\cite[Lemma 3.1]{bohm2021variable}, for some $\Lipmor > 0$.
   In particular, it holds
   \begin{align*}
      \nabla \LyapS(\x) = \tfrac{1}{\stepx} (\x - {\argmin_{\y \in \cC}} (\ell(\y) + \tfrac{1}{2\stepx} \norm{\y - \x}^2)).
   \end{align*}
   Let $\hx(\x) \coloneq \argmin_{\y \in \cC} (\ell(\y) + \tfrac{1}{2\stepx} \norm{\y - \x}^2)$.
   To lighten the notation, in the following we drop the dependence of $\hx(\x)$ and $\bax(\x)$ on $\x$.
   We proceed by studying the increment of $\LyapS$ along the trajectories of system~\eqref{eq:reduced} to verify~\eqref{eq:lyapunov_reduced_decrease}, namely
   \begin{align}
      \!\!\Delta \LyapS\left(\x\right) &\!\coloneq \LyapS\left(\x - \conv (\x - \bax)\right) - \LyapS\left(\x\right)
      \notag\\
      &\stackrel{(a)}{\leq} - \tfrac{\conv}{\stepx} \left(\x - \hx\right)\T \!\! \left(\x - \bax\right) + \tfrac{\Lipmor \conv^2}{2} \norm{\x - \bax}^2\!\!, \label{eq:delta_W}
   \end{align}
   where in $(a)$ we use the Descent Lemma~\cite[Proposition 6.1.2]{bertsekas2015convex} and the definition of $\hx$.
   By the law of cosines, we note that
   \begin{align}\label{eq:x_hx_minus_x_tx}
      2\left(\x \!-\! \hx\right)\T\!\left(\x \!-\! \bax\right) &= \norm{\x \!-\! \hx}^2 \!+\! \norm{\x \!-\! \bax}^2 \!-\! \norm{\bax \!-\! \hx}^2.
   \end{align}
   Now, by definition of $\bax$ and $\hx$ and denoting with $\Nc{\cC}[\cdot]$ the normal cone operator of the set $\cC$, we note that
   \begin{align*}
      \frac{\left(\x \!-\! \bax\right)}{\stepx} \!-\! \nabla \ell(\x) \!\in\! \Nc{\cC}[\bax]
      \quad \text{and} \quad
      \frac{\left(\x \!-\! \hx\right)}{\stepx} \!-\! \nabla \ell(\hx) \!\in\! \Nc{\cC}[\hx]. 
   \end{align*}   
   Since the set $\cC$ is closed an convex (cf. Assumption~\ref{ass:upper_cost}), the operator $\Nc{\cC}[\cdot]$ is monotone and, thus, it holds
   \begin{align*}
      \left(\tfrac{\left(\x - \bax\right)}{\stepx}  - \nabla \ell(\x) - \tfrac{\left(\x - \hx\right)}{\stepx} + \nabla \ell(\hx)\right)\T (\bax - \hx) \geq 0.
   \end{align*}   
   Then, by rearranging the terms in the above inequality and by leveraging the Lipschitz continuity of $\nabla \ell(\x)$ we get
   \begin{align*}
      -\norm{\bax - \hx}^2 
      \geq -\stepx\Lipf\norm{\bax - \hx}\norm{\x - \hx}.
   \end{align*}
   Without loss of generality, let $\bax \neq \hx$. 
   We conclude that 
   \begin{align}\label{eq:bound_tx_th}
      \norm{\bax - \hx} \leq \stepx\Lipf \norm{\x - \hx}.
   \end{align}
   By plugging~\eqref{eq:bound_tx_th} into~\eqref{eq:x_hx_minus_x_tx} and since $\stepx \in (0, \bstepxt)$, we get
   \begin{align}
      2\left(\x \!-\! \hx\right)\T\!\left(\x \!-\! \bax\right) 
      &\geq%
      \norm{\x \!-\! \bax}^2.
      \label{bound_x_hx_x_tx}
   \end{align}
   Now, we plug~\eqref{bound_x_hx_x_tx} into~\eqref{eq:delta_W} and use $\delta \leq 1$ to get%
   \begin{align*}
      \Delta \LyapS\left(\x\right) 
      \leq - \conv (\tfrac{1}{2\stepx} - \tfrac{\Lipmor}{2}) \norm{\x - \bax}^2.
   \end{align*}
   Then, for any $\ci\!>\!0$, we pick $\stepx \in (0, \bstepx)$, with $\bstepx \coloneq \min\{\bstepxt, \tfrac{1}{\Lipmor}\}$.
   Then, we satisfy~\eqref{eq:lyapunov_reduced_decrease} for all $\conv \in (0, 1)$. 
   Now, we proceed by verifying~\eqref{eq:lyapunov_reduced_nabla_bound}.
   We note that
   \begin{align*}
      \norm{\x - \hx} 
      &\!\stackrel{(a)}{\leq}\! \norm{\x - \bax} + \norm{\bax - \hx}
      \!\stackrel{(b)}{\leq}\! \norm{\x - \bax} + \stepx \Lipf \norm{\x - \hx}\!,
   \end{align*}
   where in $(a)$ we add and subtract $\bax$ and apply the triangle inequality, while in $(b)$ we used Equation~\eqref{eq:bound_tx_th}.
   It follows that
   \begin{align*}
       \norm{\x - \hx} \leq \tfrac{1}{1 - \stepx \Lipf} \norm{\x - \bax}.
   \end{align*}
   Hence, we verify~\eqref{eq:lyapunov_reduced_nabla_bound} with $\cc \coloneq \tfrac{1}{\stepx(1 - \stepx \Lipf)}$. 
   Finally,~\eqref{eq:lyapunov_reduced_lip} holds by Lipschitz continuity of $\nabla \LyapS(\x)$.

\bibliographystyle{IEEEtran}
\bibliography{bibliography}

\end{document}